\documentclass{scrartcl}
\usepackage[numbers]{natbib}
\usepackage{fontenc}
\usepackage{shadethm}
\usepackage{amsthm}
\usepackage{float}
\usepackage{bm}
\usepackage{stmaryrd}
\usepackage{mathrsfs} 
\usepackage{amsmath}
\usepackage{amssymb}
\usepackage{wasysym}
\usepackage{sectsty}
\usepackage{hyperref}
\usepackage[a4paper, left=2.7cm, right=2.7cm, top=2.5cm]{geometry}
\usepackage{chngcntr}
\counterwithin*{equation}{section}

\usepackage{cleveref}
\DeclareMathAlphabet{\mathpzc}{OT1}{pzc}{m}{it}

\sectionfont{\normalsize\centering}
\subsectionfont{\footnotesize\centering}
\newcommand{\norm}[1]{\left\lVert#1\right\rVert}

\theoremstyle{plain}
\newtheorem{thm}{Theorem}[section] % reset theorem numbering for each chapter

\theoremstyle{definition}
\newtheorem{defn}[thm]{Definition} % definition numbers are dependent on theorem numbers
\newtheorem{lem}[thm]{Lemma}
\newtheorem{prop}[thm]{Proposition}
\newtheorem{rem}[thm]{Remark}
\newtheorem{cor}[thm]{Corollary}

\def\XXint#1#2#3{{\setbox0=\hbox{$#1{#2#3}{\int}$ }
		\vcenter{\hbox{$#2#3$ }}\kern-.6\wd0}}

\usepackage[utf8]{inputenc}
\usepackage[german,english,russian]{babel}

\newcounter{MPequ}

\begin{document}\selectlanguage{english}
\begin{center}
\normalsize \textbf{\textsf{Typical field lines of Beltrami flows and boundary field line behaviour of Beltrami flows on simply connected, compact, smooth manifolds with boundary}}
\end{center}
\begin{center}
	Wadim Gerner\footnote{\textit{E-mail address:} \href{mailto:gerner@eddy.rwth-aachen.de}{gerner@eddy.rwth-aachen.de}}
\end{center}
\begin{center}
{\footnotesize	RWTH Aachen University, Lehrstuhl f\"ur Angewandte Analysis, Turmstra{\ss}e 46, D-52064 Aachen, Germany}
\end{center}
{\small \textbf{Abstract:} We characterise the boundary field line behaviour of Beltrami flows on compact, connected manifolds with vanishing first de Rham cohomology group. Namely we show that except for an at most nowhere dense subset of the boundary, on which the Beltrami field may vanish, all other field lines at the boundary are smoothly embedded $1$-manifolds diffeomorphic to $\mathbb{R}$, which approach the zero set as time goes to $\pm \infty$. 
\newline
We then drop the assumptions of compactness and vanishing de Rham cohomology and prove that for almost every point on the given manifold, the field line passing through the point is either a non-constant, periodic orbit or a non-periodic orbit which comes arbitrarily close to the starting point as time goes to $\pm \infty$. During the course of the proof we in particular show that the set of points at which a Beltrami field vanishes in the interior of the manifold is countably $1$-rectifiable in the sense of Federer and hence in particular has a Hausdorff dimension of at most $1$. As a consequence we conclude that for every eigenfield of the curl operator, corresponding to a non-zero eigenvalue, there always exists exactly one nodal domain.}
\section{Introduction}
In this paper we deal with the field line behaviour of so called Beltrami fields. Given an oriented, smooth Riemannian $3$-manifold $(\bar{M},g)$ with or without boundary, we call a smooth vector field $\bm{X}$ on $\bar{M}$ a \textbf{weak Beltrami field} if there exists a locally bounded function $\lambda:\bar{M}\rightarrow \mathbb{R}$ such that
\begin{equation}
\label{I1}
\text{curl}(\bm{X})=\lambda \bm{X}\text{ and }\text{div}(\bm{X})=0.
\end{equation}
These type of vector fields appear naturally in physics. For instance they provide stationary solutions of the equations of ideal magnetohydrodynamics for the case of constant pressure and a resting plasma, \cite[Chapter III]{AK98}, and are possible terminal magnetic field configurations for suitable initial magnetic fields. Beltrami fields, in fact, arise as local minimisers of the helicity constrained magnetic energy minimisation, \cite{A74}, \cite{AL91}, \cite{W58}. Hence Beltrami fields are of interest in ideal magnetohydrodynamics, and in astrophysics in particular. On the other hand, we recall the usual form of the incompressible Euler equations
\[
\partial_t\bm{v}+\nabla_{\bm{v}}\bm{v}=-\nabla p\quad\text{and}\quad\text{div}(\bm{v})=0,
\]
where $\bm{v}$ is the fluid velocity and $p$ is the pressure. These equations can be written equivalently as
\begin{equation}
\label{I2}
\partial_t\bm{v}+\nabla f=\bm{v}\times \text{curl}(\bm{v})\quad\text{and}\quad\text{div}(\bm{v})=0,
\end{equation}
with $f:=\frac{g(\bm{v},\bm{v})}{2}+p$ being the Bernoulli function of the system. It follows that if $\bm{v}$ is a (time-independent) weak Beltrami field, then it is a solution of (\ref{I2}) for the pressure function $p:=c-\frac{g(\bm{v},\bm{v})}{2}$ and any constant $c\in \mathbb{R}$. Therefore, weak Beltrami fields may instead be interpreted as stationary solutions of the incompressible Euler equations. In fact weak Beltrami fields are a very special kind of stationary solutions of the Euler equations. If we restrict ourselves to the case of stationary solutions of these equations, we arrive at the following PDE's
\begin{equation}
\label{I3}
\bm{v}\times \text{curl}(\bm{v})=\nabla f\text{ and }\text{div}(\bm{v})=0,
\end{equation}
for some function $f$. Assuming that all quantities involved are real analytic and that $\bm{v}$ and its curl are not everywhere collinear, \textit{Arnold's structure theorem}, \cite{A66}, \cite{A74}, \cite[Chapter II Theorem 1.2]{AK98}, \cite{GK94}, provides a rather complete picture of the corresponding fluid flow. On the other hand, if $\bm{v}$ and its curl are everywhere collinear, which includes the case of weak Beltrami fields, the topological properties of the flow can be much more complicated. For instance it is a consequence of Arnold's structure theorem, see \cite[Chapter II, proposition 6.2]{AK98}, that every no-where vanishing real analytic steady flow, which possesses a 'chaotic' trajectory, i.e. a trajectory not contained in any $2$ dimensional subset, is necessarily a Beltrami field. Hence the presence of complicated flowline topology already implies that the underlying flow must be Beltrami. Conversely there are Beltrami fields which have indeed complicated field line behaviour. The so called ABC-flows (Arnold-Beltrami-Childress) develop for certain choices of their parameters chaotic field lines, \cite{DFGHMS86}. Etnyre and Ghrist showed in \cite{EtGr00b} that there exists a metric $g$ on $S^3$ and a Beltrami field corresponding to this metric such that this Beltrami field admits all knot and link types as its fields lines. A similar result regarding the Euclidean metric was obtained by Enciso and Peralta-Salas, \cite{EP12}, who showed that there exists an eigenfield of the curl operator on $\mathbb{R}^3$ with the Euclidean metric, such that it admits all tame knots and (locally finite) links as field line configurations, see also \cite{EP15}.
These results have in common that one considers or constructs specific Beltrami fields with interesting topologies. But not all Beltrami fields have such complicated structures. Some of them behave just like 'typical' steady flows to which Arnold's structure theorem applies. For instance in \cite{CDGT00} the authors explicitly compute all Beltrami fields on the solid ball in $\mathbb{R}^3$ which are tangent to its boundary. The eigenfield corresponding to the smallest positive eigenvalue indeed admits field lines which lie on invariant tori and the ball fibres into such invariant tori, after removing a lower dimensional subset from it. This is in essence the statement of the structure theorem.
Contrary to the above results another question one may ask is: 
\newline
\newline
\textit{Is there something we can say about the field line behaviour of all Beltrami fields (or some subclass) under certain assumptions on our manifold?}
\newline
\newline
A typical question in this regards is whether or not Beltrami fields admit (non-constant) periodic orbits. Under certain assumptions the answer is positive and a consequence of the (by now proven) Weinstein conjecture, \cite{Ta07}, \cite{H093}, and the relation between Reeb vector fields and no-where vanishing, rotational Beltrami fields \cite{EtGr00a}. A rotational Beltrami field is a weak Beltrami field with smooth proportionality function $\lambda$ which is no-where zero. Weinstein's conjecture can then be stated as follows \cite{EtGr00a}: Every no-where vanishing, rotational Beltrami field on an oriented, closed (compact, without boundary) smooth, Riemannian $3$-manifold $(M,g)$ admits a periodic orbit. See also \cite{EtGr99} for results concerning the existence of closed orbits of real analytic steady flows of the Euler equations on $S^3$. Much less seems to be known about Beltrami flow behaviour on manifolds with non-empty boundary. For example \cite{EtGr01} and \cite{EtGr02} deal with the existence of closed orbits of non-singular Beltrami fields defined on and tangent to the boundary of the solid torus.
\newline
In the present paper we want to deal with the question posed above with a focus on manifolds with boundary. Our approaches are twofold.
\newline
\newline
\textit{First}: We restrict our attention to Beltrami fields tangent to the (non-empty) boundary and analyse the field line behaviour on the boundary. The idea being that it is easier to understand dynamics on surfaces rather than in $3$-space. The main observation is that the restriction of a given Beltrami field to the boundary satisfies the integrability conditions. In particular, if the boundary has a vanishing first de Rham cohomology group, the restricted vector field must be a gradient field whose dynamics are well understood. To ensure that the boundary is well-behaved we make the assumptions of compactness and $H^1_{dR}(\bar{M})=\{0\}$ on the underlying manifold $\bar{M}$. But our results also hold true on non-compact manifolds with arbitrary de Rham cohomology groups, so as long as the boundary is well-behaved.
\newline
\newline
\textit{Second:} We notice that the boundary is a $2$-dimensional subset and hence a null set. Thus the typical field lines of Beltrami fields on manifolds with and without boundary should be the same. Hence our idea is that it should be possible to characterise typical field lines. Indeed it is possible and in fact we will see that the boundary field line behaviour from the previous paragraph turns out to be atypical.
\newline
\newline
While characterising typical field lines we will show that, given a weak Beltrami field $\bm{X}$, the set of points at which it vanishes on $\bar{M}$ has Hausdorff dimension at most $2$. In addition we will see that if we restrict ourselves to the interior of the manifold and assume that the proportionality function $\lambda$ has bounded derivatives near the zero set, see \cref{MDadd} for a precise statement, then the Hausdorff dimension is at most $1$. This in particular includes the case of smooth proportionality function. The upper bound is sharp, as can be seen by considering the ABC-flows,\cite{A74}, with parameters $A=-C=1$ and $B=0$. The result regarding the Hausdorff dimension in the smooth proportionality function case is not new. In fact it is a special case of a result by B\"ar concerning general elliptic operators of first order, \cite{B99}. However B\"ar's result, even though more general, requires deeper technical tools, such as Malgrange's preparation theorem, see for instance \cite{Ni71}. Our proof only uses the implicit function theorem and a unique continuation result by exploiting the specific structure of our equations. This provides an alternative, more elementary proof for our situation in the smooth proportionality function setting. But our result allows a more general class of proportionality functions so that, in this sense, we provide a generalisation of B\"ar's results for our specific equations.
\newline
Note that if the proportionality function $\lambda$ is a constant, then the corresponding weak Beltrami field $\bm{X}$ satisfies
\[
\Delta \bm{X}=\lambda^2 \bm{X}
\]
since $\text{div}(\bm{X})=0$ and where we choose $\Delta$ to be a positive operator, i.e. it coincides in the Euclidean case with the negative of the standard Euclidean Laplacian. In the Euclidean case each component of $\bm{X}$ is therefore an eigenfunction of the Laplacian and the zero set of $\bm{X}$ is exactly the intersection of the zero sets of the component functions. Zero sets $\mathcal{N}_f$ of eigenfunctions $f$ of the Laplacian, also called nodal sets, as well as the connected components of their complement, so called nodal domains, have been widely studied in mathematics, see for instance \cite{HarSi89}, \cite{Log18}, \cite{SoZel11}, \cite{SoZel12} and see also \cite{GrbNg13} for a survey of properties of Laplacian eigenfunctions. The classical Courant's nodal domain theorem, \cite[\S 5, p.365]{CouHi24} states that the $k$-th-Dirichlet eigenfunction of the Laplacian on a bounded domain with smooth boundary in $\mathbb{R}^3$ admits at most $k$ nodal domains (it is easy to show that there must be at least two, provided $k\geq 2$). It is then natural to ask the same question with regards to the curl operator.
\newline
\newline
\textit{Suppose $\bm{X}$ is an eigenfield of the curl operator, tangent to the boundary of $\bar{M}$. How many connected components does the set $\text{int}(\bar{M})\setminus \{\bm{X}=0\}$ have?}
\newline
\newline
We will show that the answer is simple, provided the corresponding eigenvalue is non-zero: It is always one component and in fact this holds regardless of the boundary conditions we impose.
\section{Main results}
\textbf{Conventions:} Manifolds are assumed to be Hausdorff and second countable. All manifolds in question are assumed to be oriented, connected, $C^{\infty}$-smooth, with (possibly empty) boundary and are assumed to be equipped with a smooth Riemannian metric, unless specified otherwise. We will shortly write: 'Let $(\bar{M},g)$ be a $3$-manifold' and mean by that, that $(\bar{M},g)$ has all the previously listed properties and is $3$-dimensional. If $\bar{M}$ is any smooth $n$-dimensional manifold with or without boundary with $H^1_{dR}(\bar{M})=\{0\}$, i.e. vanishing first de Rham cohomology group, then we call $\bar{M}$ a simple manifold. The set of all smooth vector fields on a given manifold $\bar{M}$ is denoted by $\mathcal{V}(\bar{M})$. Given a smooth vector field $\bm{X}\in \mathcal{V}(\bar{M})$ and a point $p\in \bar{M}$ we denote by $\gamma_p$ the field line of $\bm{X}$ starting at $p$ and by $\omega(p),\alpha(p)$ the corresponding $\omega-$ and $\alpha$-limit sets. Recall that $q\in \omega(p)$ if and only if there exists a sequence $(t_n)_n$, $t_n\rightarrow+\infty$, such that $\gamma_p(t_n)\rightarrow q$ and $\alpha(p)$ is defined accordingly with limits $t_n\rightarrow-\infty$. We call a vector field $\bm{X}$ on $(\bar{M},g)$ a weak $C^k$-Beltrami field, if it is of class $C^k$ and there exists a locally bounded function $\lambda:\bar{M}\rightarrow \mathbb{R}$, such that $\text{curl}(\bm{X})=\lambda \bm{X}$ and $\text{div}(\bm{X})=0$.
\begin{thm}
\label{MT1}
Let $(\bar{M},g)$ be a compact, simple $3$-manifold with non-empty boundary $\partial\bar{M}\neq \emptyset$ and let further $\bm{X}\in \mathcal{V}(\bar{M})$ be a weak Beltrami field. If $\bm{X}$ is tangent to the boundary and not the zero vector field, then the following holds
\begin{itemize}
\item $K:=\{p\in \partial{\bar{M}}|\bm{X}(p)=0 \}$ is a closed and no-where dense subset of $\partial\bar{M}$. Equivalently $U:=\partial\bar{M}\setminus K$ is an open and dense subset of $\partial\bar{M}$.
\item The (globally defined) field lines $\gamma_p$ of $\bm{X}$ starting at points $p\in U\subset \partial\bar{M}$ satisfy the following two properties:
\begin{itemize}
\item The images $\gamma_p(\mathbb{R})\subset \partial\bar{M}$ are smoothly embedded $1$-submanifolds diffeomorphic to $\mathbb{R}$.
\item The $\omega$- and $\alpha$-limit sets of any such field line are non-empty and subsets of $K$. Both limit sets consist either of a single point or of infinitely many points, none of which is an isolated point in $\partial\bar{M}$.
\end{itemize}
\item $\#K\geq 2N$, where $\#K$ denotes the number of elements of $K$ and $N$ denotes the number of connected components of $\partial\bar{M}$. In particular $\bm{X}$ has at least two singular points on the boundary.
\end{itemize}
\end{thm}
\begin{cor}[Boundary field line behaviour of Beltrami fields on simple $3$-manifolds]
\label{MC2}
Let $(\bar{M},g)$ be a compact, simple $3$-manifold with non-empty boundary and suppose $\bm{X}\in \mathcal{V}(\bar{M})$ is a weak Beltrami field, which is tangent to the boundary and has at most finitely many non-isolated zeros in $\partial\bar{M}$. If we let $K$ and $U$ be defined as in \cref{MT1}, then for any $p\in U$, the field line $\gamma_p$ of $\bm{X}$ starting at $p$ is a smooth embedding of $\mathbb{R}$ into $\partial\bar{M}$ (and hence also into $\bar{M}$) and there exist elements $p_{\pm}\in K$ such that $\lim_{t\rightarrow \pm \infty}\gamma_p(t)=p_{\pm}$.
\end{cor}
Here we mean by (non-)isolated in $\partial\bar{M}$ that we consider these notions in the subspace topology of $\partial\bar{M}$. The fact that the set $K$ has no interior points is a consequence of the following result
\begin{lem}[Vainshtein's lemma for abstract manifolds]
\label{ML3}
Let $(\bar{M},g)$ be a $3$-manifold with non-empty boundary and $\bm{X}$ a weak $C^1$-Beltrami field. If there exists an open subset $U\subseteq \bar{M}$, such that $\partial\bar{M}\cap U\neq \emptyset$ and $\bm{X}|_{\partial\bar{M}\cap U}\equiv 0$, then $\bm{X}$ is the zero vector field, i.e. $\bm{X}\equiv 0$ on all of $\bar{M}$.
\end{lem}
A proof of \cref{ML3} for the case of bounded domains with smooth boundary in $\mathbb{R}^3$ can be found in \cite{CDGT00b}, see also \cite{V92} for the original reference. However their proof relies on the position vector field. Thus it does not immediately generalise to abstract manifolds. Here we use a different approach via the double of a manifold and a unique continuation result, inspired by the proof of \cite[Theorem 3.4.4]{S95}. In particular, we do not need to assume that $\bar{M}$ is compact.
\newline
In order to show that the boundary is well-behaved in a suitable sense, given our assumptions on $\bar{M}$, we will first prove the following result
\begin{lem}
\label{ML4}
Let $n\geq 2$ and $\bar{M}$ be a compact, orientable, smooth $n$-dimensional manifold with non-empty boundary. Then for every fixed $0\leq k \leq n-2$ we have the following implication
\begin{equation}
\label{M1}
H^{(n-1)-k}_{dR}(\bar{M})=\{0\} \Rightarrow \dim\left(H^k_{dR}\left(\partial\bar{M}\right)\right)\leq \dim\left(H^k_{dR}\left(\bar{M}\right)\right).
\end{equation}
\end{lem}
This implies the following characterisation of the boundary of compact, simple manifolds
\begin{cor}
\label{MC5}
Let $\bar{M}$ be a compact, orientable, simple, smooth $3$-dimensional manifold with non-empty boundary. Then every connected component of $\partial\bar{M}$ is diffeomorphic to $S^2$.
\end{cor}
Now let us turn to the characterisation of typical field lines of weak Beltrami flows
\begin{thm}
\label{MT6}
Let $(\bar{M},g)$ be a $3$-manifold of finite volume and let $\bm{X}\in \mathcal{V}(\bar{M})$ be a weak Beltrami field, which is tangent to the boundary of $\bar{M}$ and generates a global flow. If $\bm{X}$ is not the zero vector field, then the following holds
\begin{itemize}
\item $U_B:=\partial\bar{M}\setminus \{p\in \partial\bar{M}| \bm{X}(p)=0 \}$ is an open and dense subset of $\partial\bar{M}$.
\item Let $\mathcal{F}:=\{p\in \bar{M}| p\in \omega(p)\cap \alpha(p)\text{ and }\gamma_p \text{ is not constant} \}$. Then $\bar{M}\setminus \mathcal{F}$ is a null set.
\end{itemize}
\end{thm}
\begin{cor}[Typical field lines of weak Beltrami fields]
\label{MC7}
Let $(\bar{M},g)$ be a $3$-manifold of finite volume and $\bm{X}\in \mathcal{V}(\bar{M})$ a weak Beltrami field, which is tangent to the boundary and generates a global flow. If $\bm{X}$ is not the zero vector field, then for almost every $p\in \bar{M}$, the field line $\gamma_p$ starting at $p$ is either a non-constant periodic orbit, or $\gamma_p$ is not periodic and there exist sequences $t^{\pm}_n\rightarrow \pm \infty$ with $\lim_{n\rightarrow \infty}\gamma_p(t^{\pm}_n)=p$.
\end{cor}
Comparing \cref{MC2} with the typical field line behaviour derived in \cref{MC7} we see that the boundary field lines are either constant or non-periodic and whenever they are not periodic there do not exist any points on the field lines themselves which are $\omega$ or $\alpha$ limit points. Hence all boundary field lines in the case of simple, compact $3$-manifolds are atypical. Since the boundary is a null set, this, of course, is no contradiction.
\newline
Before stating the next result, we give a definition. Note that it immediately follows from definition of a weak Beltrami field $\bm{X}$ that the proportionality function $\lambda$ is smooth away from the set of points at which $\bm{X}$ vanishes, i.e. $\lambda\in C^{\infty}(\bar{M}\setminus K)$ with $K=\{p\in \bar{M}|\bm{X}(p)=0\}$.
\begin{defn}
\label{MDadd}
Let $(\bar{M},g)$ be a $3$-manifold and $\bm{X}\in \mathcal{V}(\bar{M})$ be a weak Beltrami field with proportionality function $\lambda$ and let $K:=\{p\in \bar{M}|\bm{X}(p)=0\}$. We say that $\lambda$ is $C^{\infty}$-\textit{bounded near} $K$ if for every $p\in K$, there exists a chart $(\mu,U)$ around $p$, such that for every multi-index $\alpha\in \mathbb{N}^3_0$ we have $\norm{\partial^{\alpha}\left(\lambda\circ \mu^{-1}\right)}_{L^{\infty}(\mu(U\setminus K))}<\infty$, where we use the standard multi-index notation.
\end{defn}
The defined notion is chart independent in the sense that if $(\mu,U)$ is any chart around some $p\in K$ with the described property and $(\tilde{\mu},\tilde{U})$ is another chart around $p$, then after shrinking $\tilde{U}$, if necessary, the $L^{\infty}$ norms of the $\tilde{\mu}$ local expressions stay finite as well. Of course, whenever $\lambda \in C^{\infty}(\bar{M})$, then it is also $C^{\infty}$-bounded near $K$. The set $U\setminus K$ may (a priori) be empty, which is for instance the case when $\lambda$ is the zero function. In that case we set the $L^{\infty}$-norm to zero.
\newline
\newline
The results of \cref{MT6} rely on the following
\begin{prop}
\label{MP8}
Let $(\bar{M},g)$ be a $3$-manifold and $\bm{X}\in \mathcal{V}(\bar{M})$ be a weak Beltrami field, which is not the zero vector field. Then the set $K:=\{p\in \bar{M}| \bm{X}(p)=0 \}$ is countably $2$-rectifiable and hence has a Hausdorff dimension of at most $2$. In particular it is a null set. If the proportionality function $\lambda:\bar{M}\rightarrow \mathbb{R}$ is $C^{\infty}$-bounded near $K$, then the set $K_I:=\{p\in \text{int}(\bar{M})| \bm{X}(p)=0 \}$ is countably $1$-rectifiable and thus has a Hausdorff dimension of at most $1$.
\end{prop}
Here we equip $(\bar{M},g)$ with the standard geodesic metric $d_g$ induced by $g$, \cite[Theorem 2.55]{L18} and the Hausdorff dimension is considered with respect to this metric. We do not require in \cref{MP8} that $\bm{X}$ is tangent to the boundary nor do we demand that $(\bar{M},g)$ has finite volume. As a consequence we obtain
\begin{cor}
\label{MC9}
Let $(M,g)$ be a $3$-manifold with empty boundary and let $\bm{X}\in \mathcal{V}(M)$ be a weak Beltrami field whose proportionality function is $C^{\infty}$-bounded near $K:=\{p\in M| \bm{X}(p)=0\}$. If $\bm{X}$ is not the zero vector field, then $M\setminus K$ is path-connected, i.e. there is exactly one nodal domain.
\end{cor}
The preceding corollary follows immediately from \cref{MP8} and the following result
\begin{prop}
\label{MP10}
Let $(M,g)$ be a connected, smooth Riemannian $n$-dimensional manifold without boundary, $n\geq 2$, and $A\subseteq M$ be a closed subset whose Hausdorff dimension is strictly less than $(n-1)$, then $M\setminus A$ is path-connected.
\end{prop}
\section{Proof of \cref{MT1}}
We will first prove \cref{ML3} and \cref{MC5} and then explain how they imply \cref{MT1}.
\newline
\newline
\underline{Proof of \cref{ML3}:} We observe that the interior of a manifold $M:=\text{int}(\bar{M})$ is homotopy equivalent to the manifold itself, \cite[Theorem 9.26]{L12}, and since connectedness is a homotopy invariant we conclude that $M$ is also connected. Now let $\hat{M}:= M\cup \left(\partial\bar{M}\cap U\right)$. Since $\partial\bar{M}\cap U$ is an open subset of the boundary, we see that $\hat{M}$ is a manifold with non-empty boundary in its own right. Observe that $M\subset \hat{M}\subseteq \bar{M}=\text{clos}(M)$, where $\text{clos}(M)$ denotes the topological closure of $M$ in $\bar{M}$. Since $M$ is connected, so must be $\hat{M}$. Now consider the boundaryless double $2\hat{M}$ of $\hat{M}$, \cite[Theorem 9.29, Example 9.32]{L12}, which is an orientable connected manifold without boundary into which $\hat{M}$ can be embedded such that the embedding is orientation preserving and such that $\hat{M}$ is a closed subset of the double. In particular the topological boundary of $\hat{M}$ in $2\hat{M}$ coincides with its manifold boundary. We can then equip the embedded copy of $\hat{M}$ with the pullback metric and extend it to a smooth metric on the whole double (the upcoming arguments are independent of the chosen extension). Then we can pushforward the vector field $\bm{X}$ to its embedded copy and extend it by $0$ outside the copy. Observe that $\text{div}(\bm{X})=0$ and hence the divergence in particular vanishes at the boundary. Further $\bm{X}$ vanishes on the boundary of $\hat{M}$ and since it is a weak Beltrami field by assumption, we see that its curl also vanishes on the boundary of $\hat{M}$. From this it follows that the vector field on the double which we extended by $0$, defines a $C^1$-vector field on the double, see also \cite[Theorem 3.4.4]{S95} for more details. One then readily checks that this extended vector field (after being identified with a $1$-form via the Riemannian metric) satisfies the requirements of a unique continuation result, \cite[Theorem 3.4.3]{S95} and \cite{AKS62}. Since it vanishes identically on an open subset of the double (the complement of $\hat{M}$) we conclude that $\bm{X}$ vanishes identically on all of $\hat{M}$ and in particular on $M$. Since $M$ is dense in $\bar{M}$ and $\bm{X}$ is continuous we see that $\bm{X}$ vanishes identically on all of $\bar{M}$. $\square$
\newline
\newline
\underline{Proof of \cref{ML4}:} For given $0\leq k\leq n-2$ we consider the restriction operator
\begin{equation}
\label{PT11}
j_k:H^k_{dR}(\bar{M})\rightarrow H^k_{dR}(\partial\bar{M}),\text{ }[\omega]\mapsto [\iota^{\#}\omega],
\end{equation}
where $\iota:\partial\bar{M}\rightarrow \bar{M}$ denotes the inclusion map and $\iota^{\#}$ its pullback. One easily checks that (\ref{PT11}) gives rise to a well-defined, linear map between the real vector spaces $H^k_{dR}(\bar{M})$ and $H^k_{dR}(\partial\bar{M})$. We claim that this operator is surjective if $H^{(n-1)-k}_{dR}(\bar{M})=\{0\}$. This will immediately imply the claim. To this end let $[\alpha]\in H^k_{dR}(\partial\bar{M})$ be any fixed element and let $\alpha\in \Omega^k(\partial\bar{M})$ be any fixed closed (smooth) $k$-form representing $[\alpha]$. By means of a collar neighbourhood construction we can find a (not necessarily closed) $k$-form $\tilde{\alpha}\in \Omega^k(\bar{M})$ such that $\iota^{\#}\tilde{\alpha}=\alpha$. We claim that there exists a \textit{closed} $k$-form $\omega\in \Omega^k(\bar{M})$ such that $\iota^{\#}\omega=\alpha$. Hence we consider the following boundary value problem
\begin{equation}
\label{PT12}
d\omega=0\text{ in }\bar{M}\text{ and }\omega|_{\partial\bar{M}}=\tilde{\alpha}|_{\partial\bar{M}}.
\end{equation}
According to \cite[Theorem 3.3.3]{S95} it has a solution if and only if the following two conditions are satisfied
\begin{equation}
\label{PT13}
0=\iota^{\#}d\tilde{\alpha}\text{ and }0=\int_{\partial\bar{M}}\iota^{\#}\left(\tilde{\alpha}\wedge \star\lambda\right),
\end{equation}
for all $\lambda\in \mathcal{H}^{k+1}_D(\bar{M})=\{\gamma\in \Omega^{k+1}(\bar{M})| d\gamma=0=\delta \gamma\text{ and }\iota^{\#}\gamma=0 \}$. Here $\delta$ denotes the adjoint derivative and $\star$ denotes the Hodge star operator. For the first condition of (\ref{PT13}) we check $\iota^{\#}d\tilde{\alpha}=d\iota^{\#}\tilde{\alpha}=d\alpha=0$, where we used that the pullback commutes with exterior derivatives, that $\iota^{\#}\tilde{\alpha}=\alpha$ and that $\alpha$ is a closed form. As for the second condition of (\ref{PT13}) we observe that $\mathcal{H}^{k+1}_D(\bar{M})\cong H^{n-(k+1)}_{dR}(\bar{M})=\{0\}$, \cite[Theorem 2.6.1, Corollary 2.6.2]{S95}, the latter equality by assumption. Therefore the second condition of (\ref{PT13}) (under our assumptions) is trivially satisfied. Hence (\ref{PT12}) admits a solution $\omega\in \Omega^k(\bar{M})$, i.e. $d\omega=0$ and $\iota^{\#}\omega=\iota^{\#}\tilde{\alpha}=\alpha$. We conclude that $\omega$ induces an equivalence class $[\omega]\in H^k_{dR}(\bar{M})$ which satisfies $j_k([\omega])=[\iota^{\#}\omega]=[\alpha]$. Thus $j_k$ is surjective. $\square$
\newline
\newline
\underline{Proof of \cref{MC5}:} Choose $n=3$, $k=1$ in \cref{ML4}, then $(n-1)-k=1$ and by assumption $\bar{M}$ is simple, i.e. $H^1_{dR}(\bar{M})=\{0\}$. Thus we conclude $\dim\left(H^1_{dR}(\partial\bar{M})\right)\leq \dim\left(H^1_{dR}(\bar{M})\right)=0$ or in other words $H^1_{dR}\left(\partial\bar{M}\right)=\{0\}$. Observe that since $\bar{M}$ is compact and orientable, so is $\partial\bar{M}$ and that $\partial\bar{M}$ is a $2$-dimensional smooth manifold without boundary. Hence any connected component $S$ of $\partial\bar{M}$ is a connected, orientable, closed, $2$-dimensional smooth manifold. By the classification theorem for compact surfaces, \cite[Chapter 9, Theorem 3.5]{H76}, we know that its diffeomorphy type is entirely characterised by its Euler characteristic. The Euler characteristic of $S$ may be computed via
\[
\chi(S)=\sum_{i=0}^2(-1)^i\dim\left(H^i_{dR}(S)\right).
\]
By connectedness of $S$ we have $\dim\left(H^0_{dR}(S)\right)=1$. Since $H^1_{dR}(\partial\bar{M})=\{0\}$, we have in particular $\dim\left(H^1_{dR}(S)\right)=0$ and finally since $S$ is compact, connected and orientable, $\dim\left(H^2_{dR}(S)\right)=1$, \cite[Theorem 17.31]{L12}. Hence we have $\chi(S)=2$ and therefore $S$ is diffeomorphic to $S^2$. $\square$
\newline
\newline
\underline{Proof of \cref{MT1}:} The first statement of \cref{MT1} is just a direct application of \cref{ML3}. The key in proving the last two statements is the following simple observation: Let $\omega^1_{\bm{X}}$ denote the $1$-form identified with $\bm{X}$ via the Riemannian metric. The fact that $\bm{X}$ is a weak Beltrami field translates to the language of forms as
\begin{equation}
\label{PT14}
\star d\omega^1_{\bm{X}}=\lambda \omega^1_{\bm{X}}\text{ and }\delta \omega^1_{\bm{X}}=0.
\end{equation}
Since we assume $\bm{X}$ to be tangent to the boundary, its normal part vanishes, i.e. in particular $n(\star d\omega^1_{\bm{X}})=\lambda n(\omega^1_{\bm{X}})=0$. Due to the duality relations $\star t=n\star $ between the normal part $n$ and the tangent part $t$, \cite[Proposition 1.2.6]{S95}, we conclude $t(d\omega^1_{\bm{X}})=0$ or equivalently $\iota^{\#}d\omega^1_{\bm{X}}=0$. Since the pullback commutes with the exterior derivative we obtain
\begin{equation}
\label{PT15}
d\iota^{\#}\omega^1_{\bm{X}}=0.
\end{equation}
Now observe that since $\bm{X}$ is tangent to the boundary, we can restrict it to the boundary and obtain a smooth vector field $\bm{X}_B\in \mathcal{V}(\partial\bar{M})$. If we now equip $\partial\bar{M}$ with the pullback metric $g_B:=\iota^{\#}g$, then we have the identity 
\[
\iota^{\#}\omega^1_{\bm{X}}=\omega^{1,g_B}_{\bm{X}_B},
\]
where $\omega^{1,g_B}_{\bm{X}_B}$ is the one form identified with $\bm{X}_B$ via the metric $g_B$ on $\partial\bar{M}$. Hence we obtain from (\ref{PT15})
\begin{equation}
\label{PT16}
d\omega^{1,g_B}_{\bm{X}_B}=0.
\end{equation}
This means that $\bm{X}_B$ satisfies the necessary integrability conditions on $\partial\bar{M}$.
\newline
\newline
\textit{Remark:} Observe that up to this point we haven't made use of the compactness and simplicity assumptions on $\bar{M}$.
\newline
\newline
Since we assume $\bar{M}$ to be compact and simple, we conclude from \cref{MC5} and (\ref{PT16}) that $\omega^{1,g_B}_{\bm{X}_B}$ is exact, i.e. there exists a smooth function $f:\partial\bar{M}\rightarrow \mathbb{R}$ such that
\begin{equation}
\label{PT17}
X_B=\text{grad}(f)\text{ with respect to the metric }g_B.
\end{equation}
The field line behaviour of gradient flows is well-understood, see for instance \cite[Part I, Chapter 1.6]{HasKa95}, so that the last two statements of \cref{MT1} follow from the gradient structure in combination with the following elementary lemma (whose proof we omit)
\begin{lem}
\label{PT1L1}
Let $M$ be a smooth $n$-dimensional manifold without boundary and let $\bm{Y}\in \mathcal{V}(M)$. Given $p\in M$, let $\gamma_p:I\rightarrow M$ denote the maximal integral curve starting at $p$ with $I=(a,b)$ for suitable $-\infty\leq a<0<b\leq +\infty$. Then $\gamma_p$ is exactly one of the following four types of curves
\begin{itemize}
\item The image $\gamma_p(I)$ is a point. In this case $I=\mathbb{R}$ and $\gamma_p$ is constant.
\item The image $\gamma_p(I)$ is a smoothly embedded, compact $1$-submanifold diffeomorphic to $S^1$. In this case $I=\mathbb{R}$ and there exists a smallest positive number $T>0$ such that $\gamma_p(t+T)=\gamma_p(t)$ for all $t\in \mathbb{R}$. $T$ is called the period of $\gamma_p$.
\item The image $\gamma_p(I)$ is a smoothly embedded $1$-submanifold diffeomorphic to $\mathbb{R}$. In this case $\gamma_p:I\rightarrow \gamma_p(I)$ is a diffeomorphism.
\item $\gamma_p(I)$ is not a smoothly embedded $1$-submanifold and not a point. In that case there exists some $x\in \gamma_p(I)$ and a sequence $(t_k)_k\subset I$ with $t_k\uparrow b$ or $t_k\downarrow a$ and $\gamma_p(t_k)\rightarrow x$ as $k\rightarrow \infty$.
\end{itemize}
\end{lem}
Observe also that we assume $\bar{M}$ to be compact and $\bm{X}$ to be tangent. Hence $\bm{X}$ has a globally defined flow, \cite[Theorem 9.34]{L12}. $\square$
\begin{rem}
\label{PT1R2}
\begin{itemize}
\item As we have seen, the restriction of $\bm{X}$ to any given boundary component always satisfies the necessary integrability conditions. Thus as long as the boundary component has a vanishing first de Rham cohomology, the restricted vector field will still turn out to be a gradient field and similar conclusions can be drawn. For instance if we consider the complement of the open unit ball in $\mathbb{R}^3$ and a weak Beltrami field on this manifold tangent to its boundary. Then its restriction to the unit sphere shows exactly the same behaviour as described in \cref{MT1}.
\item For the proof of \cref{MT1} it was enough to have \cref{ML4} at hand to conclude that $H^1_{dR}(\partial\bar{M})=\{0\}$ and hence the restricted vector field is a gradient field. Nonetheless the implication, \cref{MC5}, is worth stating.
\item \cref{MC2} is an immediate consequence of \cref{MT1}.
\end{itemize}
\end{rem}
\section{Proof of \cref{MT6}}
The first statement of \cref{MT6} is a direct consequence of \cref{ML3}. Either $\partial\bar{M}$ is empty, in which case $U_B$ is empty as well and then obviously satisfies the conclusion. Or $\partial\bar{M}$ is not empty, in which case the conclusion follows from \cref{ML3}. Thus we focus here on the second statement of \cref{MT6}. To this end let us first recall some notions for the convenience of the reader.
\begin{defn}[Null sets on manifolds]
\label{PT6D1}
A subset $N\subseteq \bar{M}$ of a smooth $n$-dimensional manifold $\bar{M}$ with or without boundary is called a \textit{null set} if there exists an atlas $\mathcal{A}$ of $\bar{M}$ such that for every chart $\mu:U\rightarrow \mu(U)$ in $\mathcal{A}$ the set $\mu(U\cap N)$ is a Lebesgue null set of $\mathbb{R}^n$. We say that some property holds \textit{almost everywhere} if there exists a null set $N\subseteq \bar{M}$ such that the property holds for all $p\in \bar{M}\setminus N$.
\end{defn}
We observe that the boundary is always mapped into a $(n-1)$-dimensional hyperplane and hence a subset $N\subseteq \bar{M}$ is a null set if and only if $N\cap \text{int}(\bar{M})$ is a null set. Thus, for the time being, we will restrict the further discussion to manifolds without boundary. Given a topological space $(X,\tau)$ we denote by $\mathcal{B}(X)$ the Borel-$\sigma$-algebra on $X$, i.e. the smallest $\sigma$-algebra on $X$ which contains all open sets.
\begin{defn}[Riemannian measure]
\label{PT6D2}
Let $(M,g)$ be an oriented, smooth Riemannian manifold without boundary of dimension $n$. Then we define the \textit{Riemannian measure} $\mu$ on $(M,g)$ via
\begin{equation}
\label{PT61}
\mu:\mathcal{B}(M)\rightarrow [0,\infty], A\mapsto \int_{M} \chi_A \omega_g,
\end{equation}
where $\omega_g$ is the Riemannian volume form and $\chi_A$ denotes the characteristic function of $A$. We call $A\in \mathcal{B}(M)$ a $\mu$-null set if $\mu(A)=0$.
\end{defn}
One easily checks that any $\mu$-null set is also a null set in the sense of \cref{PT6D1}.
\newline
The proof of \cref{MT6} is divided in two steps. In the first step we will show by means of Poincar\'{e} recurrence that the complement of the set $\mathcal{F}^{\prime}:=\{p\in \bar{M}| p\in \omega(p)\cap \alpha(p) \}$ is a null set. In a second step we will show that the set of points at which the vector field vanishes is a null set (which are exactly the constant orbits).
\newline
\newline
\underline{Step 1:} \textit{Claim:} The set $\bar{M}\setminus \{p\in \bar{M}| p\in \omega(p)\cap \alpha(p) \}$ is a null set.
\newline
\newline
\underline{Proof of claim:} Let $M:=\textit{int}(\bar{M})$ and as pointed out before we may restrict to the set of points contained in the interior since the boundary is a null set either way. Consider the Riemannian measure $\mu$, which turns $\left(M,\mathcal{B}(M),\mu\right)$ into a finite measure space, since $\mu(M)=\int_M\chi_M\omega_g=\int_M\omega_g=\text{vol}(M)<+\infty$ by assumption. Further let $\phi_t:M\rightarrow M$ denote the (global) flow of $X$, which is well-defined since diffeomorphisms map interiors to interiors. Since Beltrami fields are divergence-free by definition the induced flow is volume preserving, i.e. $(\phi_t)^{\#}\omega_g=\omega_g$ for every $t\in \mathbb{R}$. We observe that for any given $A\in \mathcal{B}(M)$ and any fixed $t\in \mathbb{R}$, we have the identity $\chi_{\phi^{-1}_{t}(A)}=\chi_A\circ \phi_t=\phi_t^{\#}\chi_A$. Making use of the fact that $\phi_t$ is volume-preserving and in particular orientation preserving we conclude
\[
\mu((\phi_t)^{-1}(A))=\int_M\chi_{\phi^{-1}_{t}(A)}\omega_g=\int_M(\phi_t^{\#}\chi_A)\phi^{\#}\omega_g=\int_M\phi_t^{\#}\left(\chi_A\omega_g\right)=\int_M\chi_A\omega_g=\mu(A),
\]
i.e. the $\phi_t$ are measure preserving transformations. Because manifolds are metrisable, \cite[Theorem 2.55]{L18}, and separable (since we assume them to be second countable) we obtain from \cite[Theorem 4]{GoHe49} that for $\mu$-almost all points $p\in M$ the point $p$ is recurrent. That is, $\mu$-almost every point $p\in M$ satisfies the following property: For every neighbourhood $U$ of $p$ there exist sequences $(t^{\pm}_k)_k$ with $t^{\pm}_k\rightarrow \pm \infty$ and $\phi_{t^{\pm}_k}(p)\in U$ for all $k\in \mathbb{N}$. Since these sequences diverge, we may pick suitable subsequences (denoted the same way) with $|t^{\pm}_k|\geq k$. Now given $j\in \mathbb{N}$ consider for any such fixed $p$ the open metric ball $B_{\frac{1}{j}}(p)$ and let $\tau^{\pm}_j$ be defined as the $j$-th element of the corresponding sequence $(t^{\pm}_k)_k$ as constructed before. Then obviously $\phi_{\tau^{\pm}_j}(p)$ converges to $p$ and by construction $\pm \tau^{\pm}_j=|\tau^{\pm}_j|\geq j$ and hence in particular $p\in \omega(p)\cap \alpha(p)$. Hence the set in question is a $\mu$-null set and every $\mu$-null set is in particular a null set in the sense of \cref{PT6D1}. This proves the claim. $\square$
\newline
\newline
\underline{Step 2:} \textit{Claim:} The set $\{p\in \bar{M}| \gamma_p\text{ is constant} \}$ is a null set.
\newline
\newline
\underline{Proof of claim:} The claim will follow immediately once we have shown \cref{MP8}.
\newline
\newline
\underline{Proof of \cref{MP8}:} The basic ideas of the proof are similar to those used for instance in \cite{HarSi89}. For the terminology used throughout, we refer to \cite[Definition 3.2.14]{Fed69}. First we give the following definition
\begin{defn}[Order of a zero]
\label{PT6D3}
Let $U\subseteq \mathbb{R}^n$ be open and $f:U\rightarrow \mathbb{R}$ be a smooth function. Let $x\in U$ be a zero of $f$, i.e. $f(x)=0$, then the \textit{order} of $x$ is defined as $\inf\{ m \in \mathbb{N}| \exists \alpha\in \mathbb{N}^n_0: |\alpha|=m\text{ and }\partial^{\alpha}f(x)\neq 0\}$, where we use the multi-index notation. If no such $m$ exists we say $x$ is a zero of infinite order. We write $\Omega_f(x)$ for the order of a given zero $x$ of $f$. If $M$ is any given smooth, $n$-dimensional manifold without boundary and $\bm{X}\in \mathcal{V}(M)$, then given a zero $p\in M$ of $\bm{X}$, we define the \textit{order} $\Omega(p)$ of $p$ as follows: Pick any chart $\mu$ of $M$ around $p$ and consider the functions $X^j:=\bm{X}(\mu^j)\circ \mu^{-1}$, the local expressions of $\bm{X}$. Then we set $\Omega(p):=\min_{1\leq j\leq n}\Omega_{X^j}(\mu(p))$.
\end{defn}
The notion of order of a given zero is independent of a particular choice of chart.
Now consider $\bm{X}$ as given in \cref{MP8} and restrict it to the interior $M:=\text{int}(\bar{M})$ of $\bar{M}$. We consider the set $K_I:=\{p\in M|\bm{X}(p)=0 \}$. Since $\bm{X}$ is a weak Beltrami field, and in particular because we assume the proportionality function $\lambda$ to be locally bounded, it is a consequence of \cite[Theorem 1]{AKS62} that all zeros have finite order. Otherwise the unique continuation result in \cite{AKS62} implies that $\bm{X}$ is the zero vector field, contradicting our assumptions. Given $n\in \mathbb{N}$ we denote by $K_n\subseteq K_I$ the set of zeros in the interior of order $n$. Then we have
\begin{equation}
\label{PT62}
K_I=\sqcup_{n\in \mathbb{N}}K_n,
\end{equation}
where $\sqcup$ indicates that this union is disjoint. Now fix any such $K_n$ and fix any point $p\in K_n$. By definition we may pick any chart $(\mu,W)$ around $p$ and consider the functions $X^j$ as in \cref{PT6D3}. Without loss of generality let $\mu(p)=0$. Then there is one $X^j$ with $\Omega_{X^j}(0)=n$ and in particular we can find a multi-index $|\beta|=n-1$ such that $\nabla (\partial^{\beta}X^j)(0)\neq 0$, where $\nabla$ denotes the standard Euclidean gradient. Define $h:=\partial^{\beta}X^j$ and observe that $\mu(K_n\cap W)\subseteq h^{-1}(0)$ by the choices we made. However $h:\mu(W)\rightarrow \mathbb{R}$ satisfies, after possibly shrinking $W$, $\nabla h(x)\neq 0$ for all $x\in \mu(W)$ and therefore $h^{-1}(0)$ defines a $2$-dimensional submanifold of $\mu(W)$. Hence $K_n\cap W\subseteq \mu^{-1}(h^{-1}(0))$, where the latter is a $2$-dimensional submanifold of $W$, since $\mu$ is a diffeomorphism between $W$ and $\mu(W)$. In particular it is also a smoothly embedded submanifold of $M$. Now given $p\in K_n$ let $W_p$ denote the corresponding neighbourhood and $H_p$ denote the corresponding $2$-dimensional hypersurface containing $K_n\cap W_p$ as constructed just now. Observe that $K_n= \bigcup_{p\in K_n}K_n\cap W_p$ where $\left(K_n\cap W_p\right)_{p\in K_n}$ forms an open cover of $K_n$. We recall that we assume $\bar{M}$ to be second countable and thus so is $K_n$ as a subset. But second countable spaces are Lindel\"{o}f and hence we can find a countable subcover $\left(K_n\cap W_j\right)_{j\in \mathbb{N}}$ of $K_n$.  Then overall we have the inclusion
\begin{equation}
\label{PT63}
K_n= \bigcup_{j\in \mathbb{N}}K_n\cap W_j\subseteq \bigcup_{j\in \mathbb{N}}H_j. 
\end{equation}
If we now fix some $H_j$ we may cover it by proper slice charts of the form $\mu: V\rightarrow B_r(0)\subseteq \mathbb{R}^3$ and $\mu\left(V\cap H_j \right)=B_r(0)\cap \{x_3=0 \}$ for suitable $V\subseteq M\subseteq \bar{M}$ open, where by proper we mean that each such chart $\mu$ extends to a chart $\tilde{\mu}:\tilde{V}\rightarrow B_{2r}(0)$ on some larger neighbourhood. It follows then from \cite[Proposition 2.51, Lemma 2.53]{L18} that the maps $\mu^{-1}:(B_r(0),d_2)\rightarrow (\bar{M},d_g)$ are Lipschitz continuous, where $d_2$ is the standard Euclidean distance and $d_g$ denotes the induced Riemannian distance. Since we have chosen $\mu$ as slice charts, we can now let $\pi:\mathbb{R}^3\rightarrow\mathbb{R}^2$ denote the projection onto the first two components and define $f:\pi(B_r(0))\rightarrow V\cap H_j$, $(x,y)\mapsto \mu^{-1}(x,y,0)$, which is Lipschitz as a composition of Lipschitz functions. Further $f$ is defined on some bounded domain of $\mathbb{R}^2$ and we have $f\left(\pi(\mu(K_n\cap V\cap H_j)) \right)=K_n\cap V\cap H_j$ so that the sets $K_n\cap V\cap H_j$ are all $2$-rectifiable. We observe that due to (\ref{PT63}) we have $K_n=\bigcup_{j\in \mathbb{N}}(K_n\cap H_j)$ and that the sets $V$ are open sets. Thus we can now argue exactly as before, that there must exist a countable open subcover of $K_n\cap H_j$ of sets of the form $K_n\cap V\cap H_j$ each of which is $2$-rectifiable. Hence $K_n$ are countably $2$-rectifiable and by (\ref{PT62}) it follows that $K_I$ is countably $2$-rectifiable.
Since $K=\{p\in \bar{M}| \bm{X}(p)=0 \}= (\partial\bar{M}\cap K)\cup K_I$ and since one can argue similarly as before, using boundary charts, that $\partial\bar{M}\cap K$ is also countably $2$-rectifiable, we conclude that $K$ is countably $2$-rectifiable and hence has Hausdorff dimension of at most $2$. This proves the first part of \cref{MP8}, which in particular concludes the proof of Step 2 and therefore the proof of \cref{MT6}. Then \cref{MC7} immediately follows from \cref{MT6}.
\newline
\newline
\textit{Proof of second part of \cref{MP8}:} Assume now that the proportionality function $\lambda$ is $C^{\infty}$-bounded near $K$. We can argue identically as before up to the point where we specify the function $h$. In particular (\ref{PT62}) still holds and we fix some $K_n$ and $p\in K_n$ and choose the multi-index $|\beta|=n-1$ with $\nabla (\partial^{\beta}X^j)(0)\neq 0$ accordingly. Without loss of generality we may assume $j=1$. We define the functions $h^i:=(\partial^{\beta}X^i)$ on $\mu(W)$ for $1\leq i\leq 3$ and $\bm{h}:=(h^1,h^2,h^3)$. Since we were allowed to choose the chart $\mu$ as we wanted, we may in particular choose it to be a normal coordinate chart centred around $p$ (recall $p$ is an interior point), see \cite[Proposition 5.24]{L18} for properties of normal coordinates. We claim that
\begin{equation}
\label{PT64}
D\bm{h}(0)\text{ has rank at least two},
\end{equation}
where $D\bm{h}(0)$ is the Jacobian matrix of $\bm{h}$ at $0$. By choice the rank of this matrix is at least one, so it is enough to exclude the possibility that the rank is exactly $1$. We recall that $\bm{X}$ is a weak Beltrami field, which implies
\begin{equation}
\label{PT65}
\text{div}(\bm{X})=0\text{ and }\text{curl}(\bm{X})=\lambda \bm{X}.
\end{equation}
The first equation reads in local coordinates
\[
0=(\det{g})^{-\frac{1}{2}}\partial_i\left(X^i\sqrt{\det{g}} \right).
\]
Applying $\partial^{\beta}$ for the $\beta$ chosen above to both sides of the above equation and using the Leibniz rule we obtain
\[
0=\sum_{\alpha\leq \beta}\binom{\beta}{\alpha}\partial^{\beta-\alpha}\left((\det{g})^{-\frac{1}{2}}\right)\partial^{\alpha}\left(\partial_i\left(X^i\sqrt{\det{g}}\right) \right).
\]
Now observe that for $|\alpha|<|\beta|$ the order of the derivatives acting upon $X^i$ in the term $\partial^{\alpha}\left(\partial_i\left(X^i\sqrt{\det{g}}\right) \right)$ is of order at most $|\beta|=n-1$. Since the order of $\bm{X}$ at $p$ is $n$, this means by definition that all these terms vanish if we evaluate them at $\mu(p)=0$. So only the term for $\alpha=\beta$ is non-zero and we obtain
\[
0=(\det{g})^{-\frac{1}{2}}\partial^{\beta}\left(\partial_i\left(X^i\sqrt{\det{g}} \right) \right)=(\partial_ih^i)(0),
\]
where we used in the last step that the additional terms from the product rule vanishes because once more the derivatives in this term acting upon $X^i$ are of order at most $|\beta|$. We find
\begin{equation}
\label{PT66}
0=(\partial_ih^i)(0)\Leftrightarrow \text{div}(\bm{h})(0)=0,
\end{equation}
where the divergence is computed with respect to the Euclidean metric. The second equation from (\ref{PT65}) reads in local coordinates as
\[
\lambda X^k=\epsilon^{ijk}\frac{\partial_i\left(X^lg_{lj}\right)}{\sqrt{\det{g}}}\text{ for all }1\leq k\leq 3.
\]
For notational simplicity we let $Y^k:=\epsilon^{ijk}\frac{\partial_i\left(X^lg_{lj}\right)}{\sqrt{\det{g}}}$. By what we have shown in the first part of the proof, the set $K$ has a Hausdorff dimension of at most $2$ and therefore does not contain any interior points. Thus there exists a sequence $(x_m)_m\subset \mu(W\setminus K)$ converging to $\mu(p)=0$. Observe that $\lambda\circ \mu^{-1}$ is smooth on $\mu(W\setminus K)$ so that we may apply the Leibniz rule to get
\[
(\partial^{\beta}Y^k)(x_m)=\sum_{\alpha\leq \beta}\binom{\beta}{\alpha}\left(\partial^{\beta-\alpha}(\lambda\circ\mu^{-1})(x_m)\right)(\partial^{\alpha}X^k)(x_m).
\]
Then the triangle inequality and $x_m\in \mu(W\setminus K)$ imply
\[
|(\partial^{\beta}Y^k)(x_m)|\leq \sum_{\alpha\leq \beta}\binom{\beta}{\alpha}\left\Vert\partial^{\beta-\alpha}(\lambda\circ\mu^{-1})\right\Vert_{L^{\infty}(\mu(W\setminus K))}|(\partial^{\alpha}X^k)(x_m)|.
\]
Since we assume $\lambda$ to be $C^{\infty}$-bounded near $K$, after possibly shrinking $W$ if necessary, we see that there is some constant $C>0$, independent of $m$, with
\[
|(\partial^{\beta}Y^k)(x_m)|\leq C\sum_{\alpha\leq \beta}|(\partial^{\alpha}X^k)(x_m)|.
\]
Taking the limit on both sides and observing that $Y^k$ as well as $X^k$ are smooth, we find
\[
|(\partial^{\beta}Y^k)(0)|\leq C\sum_{\alpha\leq \beta}|(\partial^{\alpha}X^k)(0)|=0,
\]
where we used the definition of $K_n$ and $\beta$ in the last step and that $\alpha\leq \beta$ implies $|\alpha|\leq |\beta|$. Thus $(\partial^{\beta}Y^k)(0)=0$ for $1\leq k\leq 3$. Applying the Leibniz rule to the definition of $Y^k$ we can argue similarly as before, namely that all terms vanish except for the term of highest derivative acting on the $X^l$. Since we are in normal coordinates we immediately get
\[
0=\epsilon^{ijk}\delta_{lj}(\partial_i\partial^{\beta}X^l)(0)\text{ for all }1\leq k\leq 3.
\]
Writing this out explicitly and using the definition of the $h^i$ we obtain:
\begin{equation}
\label{PT67}
(\partial_jh^i)(0)=(\partial_ih^j)(0)\text{ for all }1\leq i,j\leq 3.
\end{equation}
Now assume for a contradiction, namely that $(D\bm{h})(0)$ in (\ref{PT64}) has rank $1$. Without loss of generality we may therefore assume that $\nabla h^1(0)\neq 0$ and since the rank is $1$, both other rows must be linearly dependent, i.e. there exist constants $a_1,a_2\in \mathbb{R}$ such that $(\nabla h^2)(0)=a_2(\nabla h^1)(0)$ and $(\nabla h^3)(0)=a_3(\nabla h^1)(0)$. For notational simplicity we will drop the argument from now on. Then these relations together with (\ref{PT66}) and (\ref{PT67}) imply
\[
\partial_2h^1=\partial_1h^2=a_2\partial_1h^1\text{, }\partial_3h^1=\partial_1h^3=a_3\partial_1h^1,
\]
\[
\partial_1h^1=-\partial_2h^2-\partial_3h^3=-a_2\partial_2h^1-a_3\partial_3h^1.
\]
Now multiply the latter equation by $a_3$ and use the first two relations to conclude
\[
\partial_3h^1=a_3\partial_1h^1=-a_3a_2\partial_2h^1-a^2_3\partial_3h^1=-a^2_2(a_3\partial_1h^1)-a^2_3\partial_3h^1=(-a^2_2-a^2_3)\partial_3h^1
\]
\[
\Leftrightarrow 0=(1+a^2_2+a^2_3)\partial_3h^1\Leftrightarrow \partial_3h^1=0.
\]
From this we obtain
\[
\partial_1h^1=-a_2\partial_2h^1-a_3\partial_3h^1=-a_2\partial_2h^1=-a^2_2\partial_1h^1\Leftrightarrow 0=(1+a^2_2)\partial_1h^1\Leftrightarrow \partial_1h^1=0.
\]
Finally the latter implies $\partial_2h^1=a_2\partial_1h^1=0$ and so overall we find $(\nabla h^1)(0)=0$, a contradiction. Thus indeed the rank of $(D\bm{h})(0)$ is at least two. Without loss of generality let $(\nabla h^2)(0)$ and $(\nabla h^1)(0)$ be linearly independent. Then by the implicit function theorem we conclude that the set $L:=\{h^1=0\}\cap \{h^2=0\}$, after possibly shrinking the domain $W$, is a $1$-dimensional submanifold of $\mu(W)$ and we observe that $\mu\left(K_n \cap W \right)\subseteq L$. From here on we can argue exactly as previously that each $K_n$ is contained in a countable union of $1$-dimensional (smoothly embedded) submanifolds and therefore is countably $1$-rectifiable. The relation (\ref{PT62}) concludes the proof. $\square$
\section{Proof of \cref{MC9}}
The corollary follows immediately from \cref{MP10}. For a concise introduction to Hausdorff measures and dimensions on metric spaces see for instance \cite[Chapter 6]{MorSol95}.
\newline
\newline
\underline{Proof of \cref{MP10}:} The basic idea of the first step of the proof for the Euclidean case is taken from a math.stackexchange discussion.
\newline
\newline
\underline{Step 1:} \textit{Claim:} For every $p\in M$ there exists an open neighbourhood $U$ of $p$ such that any two points $q,\tilde{q}\in U \setminus A$ can be joined by a continuous path in $U\setminus A$.
\newline
\newline
\underline{Proof of claim:} Fix any $p\in M$, then there exists some $\epsilon(p)>0$ such that every open geodesic ball centred around $p$ whose radius is smaller than $\epsilon(p)$ is geodesically convex, \cite[Theorem 6.17]{L18}. Thus let $U$ be an open geodesic ball of radius smaller than $\epsilon(p)$ and such that its closure is a closed geodesic ball, in particular its closure is contained in a coordinate chart. Then by definition of an open geodesic ball we can find a chart $\mu:U\rightarrow B_r(0)$ for suitable $r>0$ and $\mu(p)=0$.  Observe that $\mu:U\rightarrow B_r(0)$ is Lipschitz if we equip $U$ with the induced Riemannian distance function and $\mu(U)$ with the standard Euclidean distance, \cite[Proposition 2.51,Lemma 2.53]{L18} and observe further that the induced Riemannian distance function on $U$ coincides with the restriction of the Riemannian distance function of $M$ to $U\times U$ since $U$ is geodesically convex. Now let $q,\tilde{q}\in U\setminus A$ be any two fixed points. If $\tilde{q}=q$ we may join them by the constant path, which lies entirely within $U\setminus A$. So now suppose $\tilde{q}\neq q$. Define $\tilde{y}:=\mu(\tilde{q})$, $y:=\mu(q)$ and $v:=y-\tilde{y}$. Since $\tilde{q}\neq q$ we have $v\neq 0$ and hence it spans a $1$-dimensional subspace $L$ of $\mathbb{R}^n$. Let $H$ denote the Euclidean orthogonal complement of $L$ in $\mathbb{R}^n$ and let $\pi:\mathbb{R}^n\rightarrow H$ denote the orthogonal projection onto $H$. Since $U\cap A\subseteq A$, its Hausdorff dimension is at most as large as that of $A$ and since $U\cap A\subseteq U$ the Hausdorff dimension of $U\cap A$ computed with respect to the restriction of the Riemannian distance function of $M$ to $U\times U$ coincides with the Hausdorff dimension of $U\cap A$ with respect to the Riemannian distance on $M$. In addition $\mu$ is a Lipschitz map with respect to the Euclidean metric and the restriction of the Riemannian distance of $M$ to $U\times U$. Therefore the image $\mu(U\cap A)\subseteq \mathbb{R}^n$ has a Hausdorff dimension (with respect to the Euclidean distance) strictly less than $(n-1)$. Observe that $\pi$ is a linear map between finite dimensional vector spaces and hence Lipschitz continuous. Thus $\pi\left(\mu(U\cap A)\right)\subseteq H$ has a Hausdorff dimension strictly less than $(n-1)$. But $H$ is $(n-1)$-dimensional which implies that $\pi\left(\mu(U\cap A)\right)$ does not have any interior points and so in particular $H\setminus \pi\left(\mu(U\cap A)\right)$ is dense in $H$. Further $A$ is closed in $M$, thus $U\cap A$ is closed in $U$ and since $\mu$ is a homeomorphism $\mu(U\cap A)$ is closed in $B_r(0)$. Since $q,\tilde{q}\in U\setminus A$ we see that $y,\tilde{y}\in B_r(0)\setminus \mu(U\cap A)$, which is open as the complement of a closed subset. Therefore there exist open balls $B_y$ and $B_{\tilde{y}}$ around $y$ and $\tilde{y}$ respectively which are contained in $B_r(0)\setminus \mu(U\cap A)$. Observe that $\pi(v)=0$ and thus $\pi(B_{\tilde{y}})\cap \pi(B_y)\neq \emptyset$ since $0=\pi(v)=\pi(y)-\pi(\tilde{y})$ by linearity of $\pi$ and definition of $v$. Since $\pi$ is surjective and $H$, $\mathbb{R}^n$ finite dimensional, $\pi$ is an open map. We conclude that $\pi\left(B_{\tilde{y}}\right)\cap \pi(B_y)$ is a non-empty open subset of $H$ and since $H\setminus \pi\left(\mu(U\cap A)\right)$ is dense in $H$ the intersection of $\left(H\setminus \pi\left(\mu(U\cap A)\right)\right)$ and $\pi\left(B_{\tilde{y}}\right)\cap \pi(B_y)$ is non-empty. Let $z$ be any fixed element of this intersection and let $L_z:=\pi^{-1}(\{z\})=\{z+\lambda v|\lambda \in \mathbb{R}\}$. By definition of $z$ we have $L_z\cap \mu(U\cap A)=\emptyset$ and there exist $x_0\in B_{\tilde{y}}$, $y_0\in B_y$ with $\pi(x_0)=z=\pi(y_0)$, i.e. $x_0,y_0\in L_z$ or equivalently there are $\lambda_x,\lambda_y\in \mathbb{R}$ with $x_0=z+\lambda_xv$ and $y_0=z+\lambda_yv$. Assume without loss of generality that $\lambda_x\leq \lambda_y$, then the pass $l(t):=z+tv$ with $\lambda_x\leq t\leq \lambda_y$ connects $x_0$ and $y_0$ and is contained in $L_z$ which does not intersect $\mu(U\cap A)$. In addition $x_0\in B_{\tilde{y}}\subseteq B_r(0)$ and similarly $y_0\in B_r(0)$. By convexity of $B_r(0)$ and since $l(t)$ defines a straight line connecting two points within $B_r(0)$ the whole path $l([\lambda_x,\lambda_y])$ is contained in $B_r(0)$. Finally recall that $y_0\in B_y$ and $B_y$ was an open ball around $y$ contained in $B_r(0)\setminus \mu(U\cap A)$. Hence we may connect $y_0$ and $y$ by a straight line and similarly connect $x_0$ and $\tilde{y}$ by a straight line and obtain a continuous path from $y$ to $\tilde{y}$ contained in $B_r(0)\setminus \mu(U\cap A)$. Lifting this path via $\mu^{-1}$ gives us a path contained in $U\setminus A$ connecting $q$ and $\tilde{q}$. This proves the claim.
\newline
\newline
\underline{Step 2:} \textit{Claim:} Any two points $p,q\in M\setminus A$ can be joined by a continuous path contained in $M\setminus A$.
\newline
\newline\underline{Proof of claim:} By assumption $M$ is connected and connected manifolds are in particular path connected. Thus let $\gamma:[0,1]\rightarrow M$ be a path connecting $p$ and $q$. For any given $x\in \gamma([0,1])$ let $U_x$ be a neighbourhood as in the first step, i.e. any two points within $U_x\setminus A$ can be joined by a path within $U_x\setminus A$. By compactness and connectedness of $\gamma([0,1])$ we can cover $\gamma([0,1])$ by finitely many such neighbourhoods $U_i$, $1\leq i\leq N$ such that $p\in U_1$, $q\in U_N$ and $U_i\cap U_{i+1}\neq \emptyset$ for all $i=1,\dots, N-1$, see also the proof of \cite[Proposition 1.16]{L12}. Since $A$ has a Hausdorff dimension of less than $(n-1)$ and in particular less than $n$, it does not contain any interior points, hence $M\setminus A$ is dense in $M$. But the $U_i\cap U_{i+1}$ are non-empty open subsets of $M$ and thus there exist elements $p_i\in U_i\cap U_{i+1}\cap (M\setminus A)$, $1\leq i\leq N-1$. By choice of the $U_i$ we can now connect $p$ and $p_1$ by a continuous path within $U_1\setminus A\subseteq M\setminus A$. Then connect $p_i$ with $p_{i+1}$ by a path within $U_{i+1}\setminus A\subseteq M\setminus A$ and finally join $p_{N-1}$ with $q$ within $U_N\setminus A\subseteq M\setminus A$. Overall we obtain a continuous path $\tilde{\gamma}:[0,1]\rightarrow M\setminus A$ joining $p$ and $q$. $\square$
\section*{Acknowledgements}
This work has been funded by the Deutsche Forschungsgemeinschaft (DFG, German Research Foundation) – Projektnummer 320021702/GRK2326 –  Energy, Entropy, and Dissipative Dynamics (EDDy). I would like to thank Christof Melcher and Heiko von der Mosel for discussions. Further I want to thank Daniel Peralta-Salas for pointing out the results of Christian B\"ar to me.
\bibliographystyle{plain}
\bibliography{mybibfile}
\footnotesize
\end{document}